\documentclass[english,11pt]{article}
\usepackage[latin1]{inputenc}
\usepackage{amsmath,amsthm,amssymb,babel}

\newtheorem{teo}{Theorem}
\newtheorem{defin}{Definition}

\newtheorem{lemma}{Lemma}

\def\proof{{\it Proof.}\ }
\def\endproof{\hfill $\Box$\par\vskip3mm}

\def\eq#1{(\ref{#1})}
\def\neweq#1{\begin{equation}\label{#1}}
\def\endeq{\end{equation}}
\def\weak{\rightharpoonup}
\def\eps{\varepsilon}

\def\phi{\varphi}

\def\di{\displaystyle}
\def\ri{\rightarrow}
\def\RR{{\mathbb R} }

\def\NN{{\mathbb N} }
\def\rn{{\mathbb R}^{N}}
\def\r2{{\mathbb R}^{2}}

\def\hr{\hookrightarrow}
\def\uf{\underline{f}}
\def\ug{\underline{g}}
\def\of{\overline{f}}
\def\og{\overline{g}}
\def\lequiv{\Longleftrightarrow}
\begin{document}
\title{\sc Entire solutions of Schr\"odinger elliptic systems
with discontinuous nonlinearity and sign--changing potential}
\author{\sc Teodora-Liliana Dinu \\ \small Department of Mathematics, ``Fra\c tii Buze\c sti" 
College, Bd. \c Stirbei--Vod\u a No. 5, 200352 Craiova, Romania\\ \small 
E-mail: {\tt tldinu@yahoo.com}}
  \date{}
      \maketitle
	  
\bigskip
 {\bf Abstract.} We establish the existence of an entire solution for a class of stationary 
Schr\"odinger systems with subcritical discontinuous nonlinearities and lower bounded potentials 
that blow-up at infinity. The proof is based on the critical point theory in the sense of Clarke 
and we apply Chang's version of the Mountain Pass Lemma for locally Lipschitz functionals. Our 
result generalizes in a nonsmooth framework a result of Rabinowitz \cite{rabi} related to entire 
solutions of the Schr\"odinger equation.

{\bf Key words}: nonlinear elliptic system, entire solution, Lipschitz functional, Clarke 
generalized gradient.

{\bf 2000 AMS Subject Classification}: 35J50, 49J52, 58E05.

\normalsize

\section{Introduction and the main result}
The Schr\"odinger equation plays the role of Newton's laws and conservation of energy in classical 
mechanics, that is, it predicts the future behaviour of a dynamic system. The linear form of 
Schr\"odinger's equation is
$$\Delta\psi +\frac{8\pi^2m}{\hbar^2}\,\left(E(x)-V(x)\right)\psi =0\,,$$
where $\psi$ is the Schr\"odinger wave function, $m$ is the mass, $\hbar$ denotes Planck's 
constant, $E$ is the energy, and $V$ stands for the potential energy. The structure of the 
nonlinear Schr\"odinger equation is much more complicated. This equation describes various 
phenomena arising in: 
self-channelling of a high-power
ultra-short laser in matter, in the theory of Heisenberg ferromagnets and magnons,
in dissipative quantum mechanics, in condensed matter theory, in plasma physics
(e.g., the Kurihara superfluid film equation). We refer to \cite{grosse,sulem} for a modern 
overview, including applications.

Consider the model problem
\begin{equation}\label{oheq}
i\hbar\psi_t=-\frac{\hbar^2}{2m}\,\Delta\psi+V(x)\psi-\gamma |\psi|^{p-1}\psi\qquad
\mbox{in $\RR^N$ ($N\geq 2$)}\,,\end{equation}
where $p<2N/(N-2)$ if $N\geq 3$ and $p<+\infty$ if $N=2$. In the study of this equation Oh 
\cite{oh} supposed that the potential $V$ is bounded
and possesses a non-degenerate critical point at $x=0$. More precisely, it is assumed that
$V$ belongs to the class ($V_a$) (for some $a$) introduced in Kato \cite{kato}. 
Taking $\gamma>0$ and $\hbar >0$
sufficiently small and using a Lyapunov-Schmidt type reduction, Oh \cite{oh} proved the existence 
of
a standing wave solution of Problem \eq{oheq}, that is, a solution of the form 
\begin{equation}\label{oheq1}
\psi (x,t)=e^{-iEt/\hbar}u(x)\,.\end{equation}
Note that substituting the ansatz \eq{oheq1} into \eq{oheq} leads to
$$-\frac{\hbar^2}{2}\,\Delta u+\left(V(x)-E\right)u=|u|^{p-1}u\,.$$
The  change of variable $y=\hbar^{-1}x$ (and replacing $y$ by $x$) yields
$$-\Delta u+2\left(V_\hbar (x)-E\right)u=|u|^{p-1}u\qquad\mbox{in $\RR^N$}\,,$$
where $V_\hbar (x)=V(\hbar x)$.

In a celebrated paper, Rabinowitz \cite{rabi} continued the study of standing wave solutions of 
nonlinear Schr\"odinger equations. After making
a standing wave ansatz, Rabinowitz reduces the problem to that of studying the
semilinear elliptic equation
$$
-\Delta u +b(x)u=f(x,u)\qquad\mbox{in }\RR^N,
$$
under suitable conditions on $b$ and assuming that $f$ is smooth, superlinear
and subcritical.

Inspired by Rabinowitz' paper, we 
consider the following class of coupled elliptic systems in $\rn \,(N\geq 3)$: 
\begin{equation}\label{P} \left\{\begin{array}{cc}
&\di -\Delta u_{1}+a(x)u_{1}=f(x,u_{1},u_{2})\qquad\mbox{in}\ \rn \\
&\di -\Delta u_{2}+b(x)u_{2}=g(x,u_{1},u_{2})\qquad\mbox{in}\ \rn\,.
\end{array}\right. \end{equation}

We point out that coupled nonlinear Schr\"odinger systems describe some physical phenomena such as 
the propagation in birefringent optical fibers or Kerr-like photorefractive media in optics. 
Another motivation to the study of coupled Schr\"odinger systems arises from
the Hartree-Fock theory for the double condensate, that is a binary mixture of
Bose-Einstein condensates in two different hyperfine states, cf. \cite{esry}. System \eq{P}
is also important for industrial applications in fiber communications systems \cite{hase}
and all-optical switching devices \cite{islam}.

Throughout this paper we assume that $a,\,b\in L_{\rm loc}^{\infty}({\rn})$ and there exist 
$\underline{a}\, ,\ \underline{b}>0 $ 
such that $\ a(x)\geq\underline{a}\ $, $\ b(x)\geq\underline{b}\ $ a.e. in $\rn$, 
and ${\rm esslim}_{|x|\ri\infty} a(x) ={\rm esslim}_{|x|\ri\infty} b(x)=+\infty $. 
Our aim in this paper is to study the existence of solutions to the above problem in the case
when $f,\, g$ are not continuous functions. Our goal is to show how variational methods can be 
used to find existence results for
stationary  nonsmooth Schr\"odinger systems.

Throughout this paper we assume that
$f(x,\cdot,\cdot),\, g(x,\cdot,\cdot)\in L_{\rm loc}^{\infty}(\RR^{2})$. Denote:
$$\uf(x,t_{1},t_{2})=\lim_{\delta\to 0} {\rm essinf}\{ f(x,s_{1},s_{2})\, ;\,
|t_{i}-s_{i}|\leq\delta\, ;\, i=1,2\}$$
$$\of(x,t_{1},t_{2})=\lim\limits_{\delta\to 0} {\rm esssup}\{ f(x,s_{1},s_{2})\, ;\,
|t_{i}-s_{i}|\leq\delta\, ;\, i=1,2\}$$
$$\ug(x,t_{1},t_{2})=\lim\limits_{\delta\to 0} {\rm essinf}\{g(x,s_{1},s_{2})\, ;\,
|t_{i}-s_{i}|\leq\delta\, ;\, i=1,2\}$$
$$\og(x,t_{1},t_{2})=\lim\limits_{\delta\to 0} {\rm esssup}\{g(x,s_{1},s_{2})\, ;\,
|t_{i}-s_{i}|\leq\delta\, ;\, i=1,2\}\,.$$

Under these conditions we reformulate Problem \eq{P} as follows:

\begin{equation}\label{P'} \left\{\begin{array}{cc}
&\di -\Delta u_{1}+a(x)u_{1}\in 
[\uf(x,u_{1}(x),u_{2}(x)),\of(x,u_{1}(x),u_{2}(x))]\quad\mbox{a.e.}\ x\in\rn \\
&\di -\Delta u_{2}+b(x)u_{2}\in 
[\ug(x,u_{1}(x),u_{2}(x)),\og(x,u_{1}(x),u_{2}(x))]\quad\mbox{a.e.}\ x\in\rn \,.  
\end{array}\right. \end{equation}

Let $H^{1}=H(\rn,\r2)$ denote the Sobolev space of all $U=(u_{1},u_{2})\in (L^{2}(\rn))^{2}$
with weak derivatives $\di\frac{\partial u_{1}}{\partial x_{j}},\,\di\frac{\partial 
u_{2}}{\partial x_{j}}\,
(j=1,\ldots ,N)$ also in $L^{2}(\rn)$, endowed with the usual norm
$$\|U\|_{H_{1}}^{2}=\int\limits_{\rn}(|\nabla U|^{2}+|U|^{2})\, dx=
\int\limits_{\rn}(|\nabla u_{1}|^{2}+|\nabla u_{2}|^{2}+u_{1}^{2}+u_{2}^{2})\,dx\,.$$

Given the functions $a,\,b:\rn\ri\RR$ as above, define the subspace
$$E=\{U=(u_{1},u_{2})\in H^{1}\,;\,\int\limits_{\rn}(|\nabla u_{1}|^{2}+
|\nabla u_{2}|^{2}+a(x)u_{1}^{2}+b(x)u_{2}^{2})\,dx<+\infty\}\,.$$
Then the space $E$ endowed with the norm
$$\|U\|_{E}^{2}=\int\limits_{\rn}(|\nabla u_{1}|^{2}+|\nabla 
u_{2}|^{2}+a(x)u_{1}^{2}+b(x)u_{2}^{2})\,dx$$  
becomes a Hilbert space.

Since $a(x)\geq\underline{a}>0,\, b(x)\geq\underline{b}>0$, we have 
the continuous embeddings $H^{1}\hr L^{q}(\rn,\r2)$ for all $2\leq q
\leq 2^{*}=2N/(N-2).$

We assume throughout the paper that $f,g:\rn\times\r2\ri\RR$ are nontrivial measurable functions 
satisfying the following hypotheses:

\begin{eqnarray}\label{1}
\left\{\begin{array}{ll}
&\di |f(x,t)|\leq C(|t|+|t|^{p})\ \mbox{for a.e.}\ (x,t)\in\rn\times\r2 \\
&\di |g(x,t)|\leq C(|t|+|t|^{p})\ \mbox{for a.e.}\ (x,t)\in\rn\times\r2 \,,
\end{array}\right.
\end{eqnarray}
where $p<2^*$;
\begin{eqnarray}\label{2}
\left\{\begin{array}{ll}
&\di \lim\limits_{\delta\to 0}{\rm esssup} \left\{\frac{|f(x,t)|}{|t|};\, 
(x,t)\in\rn\times(-\delta,+\delta)^{2}\right\}=0\\
&\di \lim\limits_{\delta\to 0}{\rm esssup} \left\{\frac{|g(x,t)|}{|t|};\, 
(x,t)\in\rn\times(-\delta,+\delta)^{2}\right\}=0;
\end{array}\right.
\end{eqnarray}
$f$ and $g$ are chosen so that the mapping $F:\rn\times\r2\ri\RR$ defined by
 $F(x,t_{1},t_{2}):=\int_{0}^{t_{1}} f(x,\tau,t_{2})\,d\tau +
\int_{0}^{t_{2}}g(x,0,\tau)\,d\tau$ satisfies
\begin{eqnarray}\label{3}\left\{\begin{array}{ll}
&\di F(x,t_{1},t_{2})=\int\limits_{0}^{t_{2}}g(x,t_{1},\tau)\,d\tau+
\int\limits_{0}^{t_{1}}f(x,\tau,0)\,d\tau\\
&\di\mbox{and $F(x,t_{1},t_{2})=0$ if and only if $t_{1}=t_{2}=0$};\end{array}\right.
\end{eqnarray}
there exists $\mu>2$ such that for any $x\in\RR^N$
\begin{eqnarray}\label{4}
0\leq \mu F(x,t_{1},t_{2})\leq\left\{\begin{array}{llll}
&\di t_{1}\uf(x,t_{1},t_{2})+t_{2}\ug(x,t_{1},t_{2});\ t_{1},t_{2}\in [0,+\infty)\\
&\di t_{1}\uf(x,t_{1},t_{2})+t_{2}\og(x,t_{1},t_{2});\ t_{1}\in [0,+\infty),\, 
t_{2}\in(-\infty,0]\\
&\di t_{1}\of(x,t_{1},t_{2})+t_{2}\og(x,t_{1},t_{2});\ t_{1},t_{2}\in (-\infty,0]\\
&\di t_{1}\of(x,t_{1},t_{2})+t_{2}\ug(x,t_{1},t_{2});\ t_{1}\in (-\infty,0],\, 
t_{2}\in[0,+\infty)\,.
\end{array}\right.
\end{eqnarray}

\begin{defin}\label{d1}
A function $U=(u_{1},u_{2})\in E$ is called solution to the problem \eq{P'} if there exists 
a function $W=(w_{1},w_{2})\in L^{2}(\rn,\r2)$ such that
\begin{itemize}
\item[$(i)$] $\uf(x,u_{1}(x),u_{2}(x))\leq w_{1}(x)\leq\of(x,u_{1}(x),u_{2}(x))$ a.e. $x$ in 
$\rn$;\\
 $\ug(x,u_{1}(x),u_{2}(x))\leq w_{2}(x)\leq\og(x,u_{1}(x),u_{2}(x))$ a.e. $x$ in $\rn$;
\item[$(ii)$] $\di\int\limits_{\rn} (\nabla u_{1}\nabla v_{1}+\nabla u_{2}\nabla v_{2}+
a(x)u_{1}v_{1}+b(x)u_{2}v_{2})\,dx=\di\int\limits_{\rn}(w_{1}v_{1}+w_{2}v_{2})\,dx$, for all 
$(v_{1},v_{2})\in E$. 
\end{itemize}
\end{defin}

Our main result is the following.

\begin{teo}\label{t1}
Assume that conditions (\ref{1}) - (\ref{4}) are fulfilled. Then Problem \eq{P'} has at least a 
nontrivial 
solution in $E$.
\end{teo}

\section{Auxiliary results}
We first recall some basic notions from the Clarke gradient theory
for locally Lipschitz functionals (see \cite{3,4} for more details). Let $E$ be a real Banach 
space and assume that $I:E\ri\RR$ is a locally Lipschitz functional. Then the Clarke
generalized gradient is defined by
$$\partial I(u)=\{\xi\in E^{*};\, I^{0}(u,v)\geq\langle\xi,v\rangle\,,
\ \mbox{for all}\ v\in E\}\,,$$
where $I^{0}(u,v)$ stands for the directional derivative of $I$ at $u$ in the direction $v$, that 
is,
$$ I^{0}(u,v)=\limsup\limits_{w\to u\atop \lambda\searrow 0}\di\frac{I(w+
\lambda v)-I(w)}{\lambda}\,.$$

Let $\Omega$ be an arbitrary domain in $\rn$. Set
$$E_{\Omega}=\left\{U=(u_{1},u_{2})\in H^{1}(\Omega;\r2)\,;\,\int\limits_{\Omega}
(|\nabla u_{1}|^{2}+|\nabla u_{2}|^{2}+a(x)u_{1}^{2}+b(x)u_{2}^{2})\,dx<+\infty\right\}$$
which is endowed with the norm
$$\|U\|_{E_{\Omega}}^{2}=\int\limits_{\Omega}(|\nabla u_{1}|^{2}+|\nabla 
u_{2}|^{2}+a(x)u_{1}^{2}+b(x)u_{2}^{2})\,dx\,.$$  
Then $E_{\Omega}$ becomes a Hilbert space.

\begin{lemma}\label{l1}
The functional $\Psi_{\Omega}:E_{\Omega}\ri\RR, $$\Psi_{\Omega}(U)=
\int\limits_{\Omega}F(x,U)\,dx$ is locally Lipschitz on $E_{\Omega}$.
\end{lemma}

\proof. We first observe that 
$$F(x,U)=F(x,u_{1},u_{2})=\int\limits_{0}^{u_{1}} f(x,\tau,u_{2})\,d\tau +
\int\limits_{0}^{u_{2}}g(x,0,\tau)\,d\tau=\int\limits_{0}^{u_{2}}g(x,u_{1},\tau)\,d\tau+
\int\limits_{0}^{u_{1}}f(x,\tau,0)\,d\tau$$ is a Carath\'eodory functional
which is locally Lipschitz with respect to the second variable. Indeed, by
\eq{1}

$$ |F(x,t_{1},t)-F(x,s_{1},t)|=\bigg|\int\limits_{s_{1}}^{t_{1}} f(x,\tau ,t)\,d\tau\bigg|
\leq\bigg|\int\limits_{s_{1}}^{t_{1}}C(|\tau ,t|+|\tau ,t|^{p})\,d\tau\bigg|\leq
k(t_{1},s_{1},t)|t_{1}-s_{1}|\,.$$
Similarly
$$ |F(x,t,t_{2})-F(x,t,s_{2})|\leq k(t_{2},s_{2},t)|t_{2}-s_{2}|\,.
$$
Therefore
$$\begin{array}{ll}
\di |F(x,t_{1},t_{2})-F(x,s_{1},s_{2})|&\leq|F(x,t_{1},t_{2})-F(x,s_{1},t_{2})|
+|F(x,t_{1},s_{2})-F(x,s_{1},s_{2})|\\
&\di \leq k(V)|(t_{2},s_{2})-(t_{1},s_{1})|\,,\end{array}$$
where $V$ is a neighbourhood of $(t_{1},t_{2}),\,(s_{1},s_{2})$.

Set $$\chi_{1}(x)=\max\{u_{1}(x),v_{1}(x)\}\,\qquad 
\chi_{2}(x)=\max\{u_{2}(x),v_{2}(x)\}\quad\mbox{for all } x\in \Omega\,.$$ 
It is obvious that if $U=(u_{1},u_{2}),\,V=(v_{1},v_{2})$ belong to $E_{\Omega}$, then 
$(\chi_{1},\chi_{2})\in E_{\Omega}$. So, by H\"older's inequality and the continuous embedding 
 $E_{\Omega}\subset L^{p}(\Omega ;\r2)$,
$$ |\Psi_\Omega (U)-\Psi_\Omega (V)|\leq 
C(\|\chi_{1},\chi_{2}\|_{E_{\Omega}})\|U-V\|_{E_{\Omega}}\,,$$
which concludes the proof.
\endproof

The following result is a generalization of Lemma 6 in \cite{8}.

\begin{lemma}\label{l2}
Let $\Omega$ be an arbitrary domain in $\rn$ and let $f:\Omega\times\r2\ri\RR$ be
a Borel function such that $f(x,.)\in L_{\rm loc}^{\infty}(\r2)$. Then $\uf$ and $\of$
 are Borel functions.
\end{lemma}

\proof Since the requirement is local we may suppose that $f$ is bounded by $M$
and it is nonnegative. Denote by
$$f_{m,n}(x,t_{1},t_{2})=\bigg(\int\limits_{t_{1}-\frac{1}{n}}^{t_{1}+\frac{1}{n}}
\int\limits_{t_{2}-\frac{1}{n}}^{t_{2}+\frac{1}{n}}|f(x,s_{1},s_{2})|^{m}
\,ds_{1}ds_{2}\bigg)^{\frac{1}{m}}\,.$$
Since $\of(x,t_{1},t_{2})=\lim\limits_{\delta\to 0} {\rm esssup}\{ f(x,s_{1},s_{2})\, ;\,
|t_{i}-s_{i}|\leq\delta\, ;\, i=1,2\}$ we deduce that for every $\eps >0$, there exists $n\in  
\NN^{*}$ such that for 
$|t_{i}-s_{i}|<\frac{1}{n}\ (i=1,2)$ we have $|{\rm esssup} f(x,s_{1},s_{2})-\of(x,t_{1},t_{2})|
<\eps$ or, equivalently,
\begin{eqnarray}\label{5}
\of(x,t_{1},t_{2})-\eps<{\rm esssup} f(x,s_{1},s_{2})<\of(x,t_{1},t_{2})+\eps\,.
\end{eqnarray}
By the second inequality in \eq{5} we obtain
$$f(x,s_{1},s_{2})\leq\of(x,t_{1},t_{2})+\eps\quad\mbox{a.e.} \ x\in\Omega\quad
\mbox{for}\ |t_{i}-s_{i}|<\frac{1}{n}\ (i=1,2)$$
which yields 
\begin{eqnarray}\label{6}
f_{m,n}(x,t_{1},t_{2})\leq(\of(x,t_{1},t_{2})+\eps)\bigg(\sqrt{4/n^{2}}\bigg)^{\frac{1}{m}}\,.
\end{eqnarray}
Let $$A=\bigg\{ (s_{1},s_{2})\in\r2\ ; \ |t_{i}-s_{i}|<\frac{1}{n} \ (i=1,2)\, ;\
\of(x,t_{1},t_{2})-\eps\leq f(x,s_{1},s_{2})\bigg\}\,.$$
By the first inequality in \eq{5} and the definition of the essential supremum
we obtain that $|A|>0$ and
\begin{eqnarray}\label{7}
f_{m,n}\leq\bigg(\int\int\limits_{A}(f(x,s_{1},s_{2}))^{m}\, ds_{1}\, ds_{2}\bigg)^{\frac{1}{m}}
\geq(\of(x,s_{1},s_{2})-\eps)\,|A|^{1/m}\,.
\end{eqnarray}
Since \eq{6} and \eq{7} imply
$$\of(x,t_{1},t_{2})=\lim\limits_{n\to\infty}\lim\limits_{m\to\infty}f_{m,n}(x,t_{1},t_{2})\,,$$
it suffices to prove that $f_{m,n}$ is Borel.
Let
$$
{\cal M}=\{ f:\Omega\times\r2\ri\RR;\ |f|\leq M\ \mbox{and}\  f\ \mbox{is a Borel function}\}$$
$${\cal N}=\{ f\in{\cal M};\ f_{m,n}\ \mbox{is a Borel function}\}\,.
$$
Cf. \cite[p.178]{berbe}, ${\cal M}$ is the smallest set of functions having the following
properties:
\begin{itemize}
\item[$(i)$] $\{f\in C(\Omega\times\r2; \RR);\ |f|\leq M\}\subset {\cal M}$;
\item[$(ii)$] $f^{(k)}\in{\cal M}$ and $f^{(k)} \stackrel{k}{\to}f$ imply $f\in{\cal M}$.
\end{itemize}
Since ${\cal N}$ contains obviously the continuous functions and $(ii)$ is also
true for ${\cal N}$ then, by the Lebesgue dominated convergence theorem, we obtain that ${\cal 
M=N}$.
For $\uf$ we note that $\uf=-(-\of))$ and the proof of Lemma \ref{l2} is complete.
\endproof

Let us now assume that $\Omega\subset\rn$ is a bounded domain. By the continuous
embedding $L^{p+1}(\Omega ;\,\r2)\hr L^{2}(\Omega ;\r2)$, we may define the locally Lipchitz 
functional
$\Psi _{\Omega}:L^{p+1}(\Omega ;\r2)\ri\RR$ by  $\Psi _{\Omega}(U)=\di\int\limits_{\Omega} 
F(x,U)\,dx$.

\begin{lemma}\label{l3}
Under the above assumptions and for any $U\in L^{p+1}(\Omega ;\r2)$, we have
$$\partial\Psi_{\Omega}(U)(x)\subset[\uf(x,U(x)),\of(x,U(x))]\times
[\ug(x,U(x)),\og(x,U(x))]\qquad\mbox{a.e.} \ x\ \mbox{in}\ \Omega\,,$$
in the sense that if $\ W=(w_{1},w_{2})\in\partial\Psi_{\Omega}(U)\subset
L^{p+1}(\Omega ;\r2)$ then
\begin{eqnarray}
\uf(x,U(x))\leq w_{1}(x)\leq\of(x,U(x))\qquad\mbox{a.e.}\  x\ \mbox{in}\ \Omega \label{8}\\
\ug(x,U(x))\leq w_{2}(x)\leq\og(x,U(x))\qquad\mbox{a.e.}\ x\ \mbox{in}\ \Omega\,. \label{9}
\end{eqnarray}
\end{lemma}

\proof By the definition of the Clarke gradient we have
$$\int\limits_{\Omega}(w_{1}v_{1}+w_{2}v_{2})\,dx\leq\Psi_{\Omega}^{0}(U,V)\qquad
\mbox{for all}\ V=(v_{1},v_{2})\in L^{p+1}(\Omega ;\r2)\,.$$
Choose $V=(v,0)$ such that $v\in L^{p+1}(\Omega)$, $v\geq 0$ a.e. in $\Omega$. Thus, by
Lemma \ref{l2},
\begin{equation}\label{10}\begin{array}{ll}
\di \int\limits_{\Omega}w_{1}v&\di\leq \limsup\limits_{(h_{1},h_{2})\to U\atop 
\lambda\searrow 0}\frac{\int\limits_{\Omega}\bigg(\int\limits_{h_{1}(x)}
^{h_{1}(x)+\lambda v(x)}f(x,\tau,h_{2}(x)\bigg)\, dx}{\lambda}\\
&\di \leq\int\limits_{\Omega}\bigg(\limsup\limits_{(h_{1},h_{2})\to U\atop\lambda\searrow 
0}\frac{1}{\lambda}
\int\limits_{h_{1}(x)}^{h_{1}(x)+\lambda v(x)}f(x,\tau,h_{2}(x)\bigg)\, dx\leq
\int\limits_{\Omega}\of(x,u_{1}(x),u_{2}(x))v(x)\,dx\,.
\end{array}\end{equation}
Analogously we obtain
\begin{eqnarray}\label{11}
\int\limits_{\Omega}\uf(x,u_{1}(x),u_{2}(x))v(x)\,dx \leq\int\limits_{\Omega}w_{1}v\,dx
\qquad\mbox{for all}\ v\geq 0\ \mbox{in}\ \Omega .
\end{eqnarray}

Arguing by contradiction, suppose that \eq{8} is false. Then there exist $\eps >0$, a set 
$A\subset\Omega $
with $|A|>0$ and $w_{1}$ as above such that
\begin{eqnarray}\label{12}
w_{1}(x)>\of(x,U(x))+\eps\qquad\mbox{in}\ A\,.
\end{eqnarray}
Taking $v={\bf 1}_{A}$ in \eq{10} we obtain 
$$\int\limits_{\Omega}w_{1}v\,dx=\int\limits_{A}w_{1}\, dx\leq\int\limits_{A}\of(x,U(x))\, dx\,,$$
which contradicts \eq{12}. Proceeding in the same way we obtain the corresponding result for $g$ 
in \eq{9}.
\endproof

By Lemma \ref{l3}, Lemma 2.1 in Chang \cite{2} and the embedding $E_{\Omega}\hr 
L^{p+1}(\Omega,\r2)$
we obtain also that for $\Psi _{\Omega}:E_{\Omega}\ri\RR$, $\Psi 
_{\Omega}(U)=\di\int\limits_{\Omega} F(x,U)\,dx$
we have
$$\partial\Psi_{\Omega}(U)(x)\subset[\uf(x,U(x)),\of(x,U(x))]\times
[\ug(x,U(x)),\og(x,U(x))]\qquad\mbox{a.e.} \ x\in \Omega\,.$$
Let $V\in E_{\Omega}$. Then $\tilde{V}\in E$, where $\tilde{V}:\rn\ri\r2$ is defined by
$$\tilde{V}=\left\{\begin{array}{ll}
&\di V(x)\quad x\ \mbox{in}\ \Omega\\
&\di 0\quad\mbox{otherwise}\,.
\end{array}\right.$$
 For $W\in E^{*}$ we consider $W_{\Omega}\in E^{*}_{\Omega}$
such that $\langle W_{\Omega},V\rangle=\langle W,\tilde V\rangle$ for all $V$ in $E_{\Omega}$.
Set $\Psi:E\ri\RR$, $\Psi(U)=\di\int\limits_{\rn} F(x,U)$.

\begin{lemma}\label{l4}
Let $W\in\partial\Psi(U)$, where $U\in E$. Then $W_{\Omega}\in\partial\Psi_{\Omega}(U)$, 
in the sense that $W_{\Omega}\in\partial\Psi_{\Omega}(U|_{\Omega})$.
\end{lemma}

\proof By the definition of the Clarke gradient we deduce that $\langle W,\tilde{V}\rangle
\leq\Psi^{0}(U,\tilde V)$ for all $V$ in $E_{\Omega}$
$$\Psi^{0}(U,\tilde V)=\limsup\limits_{H\to U,\,H\in E\atop \lambda\to 0}
\frac{\Psi(H+\lambda\tilde{V})-\Psi(H)}{\lambda}=\limsup\limits_{H\to U,\,H\in E\atop \lambda\to 
0} 
\frac{\int\limits_{\rn}(F(x,H+\lambda\tilde{V})-F(x,H))\,dx}{\lambda}=$$
$$\limsup\limits_{H\to U,\,H\in E\atop \lambda\to 0}
\frac{\int\limits_{\Omega}(F(x,H+\lambda\tilde{V})-F(x,H))dx}{\lambda}=
\limsup\limits_{H\to U,\,H\in E_{\Omega}\atop \lambda\to 0}
\frac{\int\limits_{\Omega}(F(x,H+\lambda\tilde{V})-F(x,H))dx}{\lambda}=\Psi^{0}_{\Omega}(U,V)\,.$$ 
Hence $\langle W_{\Omega},V\rangle\leq\Psi_{\Omega}^{0}(U,V)$ which implies 
$W_{\Omega}\in\partial\Psi^{0}_{\Omega}(U)$.
\endproof

By Lemmas \ref{l3} and \ref{l4} we obtain that for any $W\in\partial\Psi(U)$ (with $U\in E$), 
$W_{\Omega}$ satisfies
\eq{8} and \eq{9}. We also observe that for $\Omega_{1},\, \Omega_{2}\subset\rn$
we have $W_{\Omega_{1}}|_{\Omega_{1}\cap\Omega_{2}}=W_{\Omega_{2}}|_{\Omega_{1}\cap\Omega_{2}}$.

Let  $W_{0}:\rn\ri\RR$, where $W_{0}(x)=W_{\Omega}(x)$ if $x\in\Omega$. Then
$W_{0}$ is well defined and
$$W_{0}(x)\in[\uf(x,U(x)),\of(x,U(x))]\times[\ug(x,U(x)),\og(x,U(x))]
\qquad\mbox{a.e.}\ x \in \rn$$
and, for all $\phi\in C_{c}^{\infty}(\rn,\r2)$,  $\langle W,\phi\rangle=\di\int\limits_{\rn} 
W_{0}\phi$.
By density of $C_{c}^{\infty}(\rn,\r2)$ in $E$ we deduce  that
$\langle W,V\rangle=\di\int\limits_{\rn}W_{0} V\,dx$ for all $V$ in $E$. Hence
\begin{eqnarray}\label{13}
W(x)=W_{0}(x)\in[\uf(x,U(x)),\of(x,U(x))]\times[\ug(x,U(x)),\og(x,U(x))]\qquad\mbox{a.e. } x \in  
\rn\,.
\end{eqnarray}

\section{Proof of Theorem \ref{t1}}
Define the energy functional $I:E\ri\RR$
\begin{eqnarray}\label{14}
I(U)=\frac{1}{2}\int\limits_{\rn}\left(|\nabla u_{1}|^{2}+|\nabla u_{2}|^{2}+
a(x)u_{1}^{2}+b(x)u_{2}^{2}\right)\,dx -\int\limits_{\rn}F(x,U)\, 
dx=\frac{1}{2}\|U\|_{E}^{2}-\Psi(U)\,.
\end{eqnarray}

The existence of solutions to problem \eq{P'} will be justified by a nonnsmooth 
variant of the Mountain-Pass Theorem (see \cite{2}) applied to the 
functional $I$, even if the $PS$ condition is not fulfilled. More precisely, we 
check the following geometric hypotheses:

\begin{eqnarray}
&\di I(0)=0\quad\mbox{and there exists} \ V\in E\quad\mbox{such that}\ I(V)\leq 0;\label{15}\\
&\di \mbox{there exist}\ \beta,\rho>0\quad\mbox{such that}\ I\geq\beta\quad\mbox{on}\ \{U\in E;\ 
\|U\|_{E}=\rho\}. \label{16}
\end{eqnarray}

{\sc Verification of \eq{15}.} It is obvious that $I(0)=0$. For the second assertion
we need the following lemma.

\begin{lemma}\label{l5}
There exist two positive constants $C_{1}$ and $C_{2}$ such that
$$f(x,s,0)\geq C_{1}s^{\mu -1}-C_{2}\quad\mbox{for a.e.}\ x\in \rn;\, s\in[0,+\infty)\,.$$
\end{lemma}

\proof We first observe that \eq{4} implies
$$0\leq\mu F(x,s,0)\leq\left\{\begin{array}{ll}
&\di s\uf(x,s,0)\qquad s\in[0,+\infty)\\
&\di s\of(x,s,0)\qquad s\in(-\infty,0]\,,
\end{array}\right.$$
which places us in the conditions of Lemma 5 in \cite{8}.

{\sc Verification of \eq{15} continued.} Choose $v\in C^{\infty}_{c}(\rn)-\{0\}$
so that $v\geq 0$ in $\rn$. We have $\int\limits_{\rn}|\nabla v|^{2}
+a(x)v^{2}<\infty$, hence $t(v,0)\in E$ for all $t\in\RR$. Thus by Lemma \ref{l5} we obtain
$$\begin{array}{ll}
\di I(t(v,0))&\di =\frac{t^{2}}{2}\int\limits_{\rn}|\nabla v|^{2}+a(x)v^{2}\,dx-\int
\limits_{\rn}\int\limits_{0}^{tv}f(x,\tau,0)\,d\tau\\
&\di\leq \frac{t^{2}}{2}\int\limits_{\rn}|\nabla v|^{2}+a(x)v^{2}\,dx-
\int\limits_{\rn}\int\limits_{0}^{tv}(C_{1}\tau^{\mu -1}-C_{2})\,d\tau\\
&\di = \frac{t^{2}}{2}\int\limits_{\rn}|\nabla v|^{2}+a(x)v^{2}\,dx+
C_{2}t\int\limits_{\rn}v\,dx-C'_{1}t^{\mu}\int\limits_{\rn}v^{\mu}\,dx<0
\end{array}$$
for $t>0$ large enough.

{\sc Verification of \eq{16}.} We observe that \eq{2}, \eq{3} and \eq{4} 
imply that, for any $\eps >0$, there exists a constant $A_{\eps}>0$ such that
\begin{eqnarray}\label{17}
\begin{array}{cc}
&\di |f(x,s)|\leq\eps|s|+A_{\eps}|s|^{p}\\
&\di |g(x,s)|\leq\eps|s|+A_{\eps}|s|^{p}
\end{array}
\qquad\mbox{for a.e. } (x,s)\in\rn\times\r2\,.
\end{eqnarray}
By \eq{17} and Sobolev's embedding theorem we have, for any $U\in E$,
$$\begin{array}{ll}
\di |\Psi(U)|&\di=|\Psi(u_{1},u_{2})|\leq\int\limits_{\rn}\int\limits_{0}^{|u_{1}|}|
|f(x,\tau,u_{2})|d\tau+\int\limits_{\rn}\int\limits_{0}^{u_{2}}|g(x,0,\tau)|\,d\tau\\
&\di\leq 
\int\limits_{\rn}\left(\frac{\eps}{2}|(u_{1},u_{2})|^{2}+\frac{A_{\eps}}{p+1}|(u_{1},u_{2}|^{p+1}
\right)\, dx+
\int\limits_{\rn}\left(\frac{\eps}{2}|u_{2}|^{2}+\frac{A_{\eps}}{p+1}|u_{2}|^{p+1}\right)\, dx\\
&\di\leq \eps\|U\|_{L^{2}}^{2}+\frac{2A_{\eps}}{p+1}\|U\|_{L^{p+1}}^{p+1}\leq
\eps C_{3}\|U\|_{E}^{2}+C_{4}\|U\|_{E}^{p+1}\,,
\end{array}$$
where $\eps$ is arbitrary and $C_{4}=C_{4}(\eps)$. Thus
$$I(U)=\frac{1}{2}\|U\|_{E}^{2}-\Psi(U)\geq\frac{1}{2}\|U\|_{E}^{2}-
\eps C_{3}\|U\|_{E}^{2}-C_{4}\|U\|_{E}^{p+1}\geq\beta> 0\,,$$
for $\|U\|_{E}=\rho$, with $\rho,\ \eps$ and $\beta$ sufficiently small positive constants.

Denote
$${\cal P}=\{\gamma\in C([0,1],E);\ \gamma(0)=0,\ \gamma(1)\not=0 \quad\mbox{and}\ 
I(\gamma(1))\leq 0\}$$
and
$$c=\inf_{\gamma\in{\cal P}}\max_{t\in[0,1]}I(\gamma(t))\,.$$
Set $$\lambda_{I}(U)=\min_{\xi\in\partial I(U)}\|\xi\|_{E^{*}}\,.$$
Thus, by the nonsmooth version of the Mountain Pass Lemma \cite{2}, there exists a sequence 
$\{U_{M}\}\subset E$ such that
\begin{eqnarray}\label{18}
I(U_{m})\to c\qquad\mbox{and}\qquad \lambda_{I}(U_{m})\to 0\,.
\end{eqnarray}
So, there exists a sequence $\{W_{m}\}\subset\partial\Psi(U_{m});\ W_{m}=(w_{m}^{1}, w_{m}^{2})$
such that
\begin{eqnarray}\label{19}
(-\Delta u_{m}^{1}+a(x)u_{m}^{1}-w_{m}^{1},-\Delta u_{m}^{2}+a(x)u_{m}^{2}-w_{m}^{2})\to 
0\quad\mbox{in}\ E^{*}\,.
\end{eqnarray}
Note that, by \eq{4},
$$\Psi(U)\leq\frac{1}{\mu}\left(\int\limits_{u_{1}\geq 0}u_{1}(x)\uf(x,U)dx+
\int\limits_{u_{1}\leq 0}u_{1}(x)\of(x,U)dx+\int\limits_{u_{2}\geq 0}u_{1}(x)\ug(x,U)dx\
+\int\limits_{u_{2}\leq 0}u_{2}(x)\og(x,U)dx\right).$$
Therefore, by \eq{13}, 
$$\Psi(U)\leq\frac{1}{\mu}\int\limits_{\rn}U(x)W(x)\, dx=\frac{1}{\mu}
\int\limits_{\rn}(u_{1}w_{1}+u_{2}w_{2})\, dx\,,$$
for every $U\in E$ and $W\in\partial\Psi(U)$.
Hence, if $\langle\cdot,\cdot\rangle$ denotes the duality pairing between $E^{*}$
and $E$, we have
$$\begin{array}{cccc}
&\di I(U_{m})=\frac{\mu-2}{2\mu}\int\limits_{\rn}(|\nabla u_{m}^{1}|^{2}+|\nabla u_{m}|^{2}+
a(x)|u_{m}|^{1}+b(x)|u_{m}|^{2})\, dx+\\
&\di \frac{1}{\mu}\langle(-\Delta u_{m}^{1}+a(x)u_{m}^{1}-w_{m}^{1},-\Delta 
u_{m}^{2}+b(x)u_{m}^{2}-w_{m}^{2}), U_{m}\rangle+\\
&\di \frac{1}{\mu}\langle W_{m},U_{m}\rangle-\Psi(U_{m})\geq\frac{\mu-2}{2\mu}
\int\limits_{\rn}(|\nabla u_{m}^{1}|^{2}+|\nabla 
u_{m}^{2}|^{2}+a(x)|u_{m}^{1}|^{2}+b(x)|u_{m}^{2}|^{2})\, dx+\\
&\di \frac{1}{\mu}\langle(-\Delta u_{m}^{1}+a(x)u_{m}^{1}-w_{m}^{1},
-\Delta u_{m}^{2}+b(x)u_{m}^{2}-w_{m}^{2}),U_{m}\rangle\geq\frac{\mu -2}{2\mu}\|U_{m}\|_{E}^{2}
-o(1)\|U_{m}\|_{E}\,.
\end{array}$$
This relation in conjunction  with \eq{18} implies that the Palais-Smale sequence $\{U_{m}\}$
is bounded in $E$. Thus, it converges weakly (up to a subsequence) in $E$ and strongly
in $L_{\rm loc}^{2}(\RR^N)$ to some $U$. Taking into account that $W_{m}\in\partial\Psi(U_{m})$
 and $U_{m}\weak U$ in $E$, we deduce from \eq{19}  that
there exists $W\in E^{*}$ such that $W_{m}\weak W$ in $E^{*}$ (up to a subsequence).
Since the mapping $U\longmapsto F(x,U)$ ia compact from  $E$ to $L^{1}$, it follows that
$W\in\partial\Psi(U)$. Therefore
$$W(x)\in[\uf(x,U(x)),\of(x,U(x))]\times[\ug(x,U(x)),\og(x,U(x))]\quad\mbox{a.e.}\ x \ \mbox{in}\ 
\rn$$
and
$$\begin{array}{cc}
&\di (-\Delta u_{m}^{1}+a(x)u_{m}^{1}-w_{m}^{1},-\Delta 
u_{m}^{2}+b(x)u_{m}^{2}-w_{m}^{2})=0\lequiv\\
&\di \int\limits_{\rn} (\nabla u_{1}\nabla v_{1}+\nabla u_{2}\nabla v_{2}+
a(x)u_{1}v_{1}+b(x)u_{2}v_{2})dx=\int\limits_{\rn}(w_{1}v_{1}+w_{2}v_{2})dx\ \mbox{for all}\ 
(v_{1},v_{2})\in E.
\end{array}$$ 
These last two relations show that $U$ is a solution pf the problem \eq{P'}.

It remains to prove that $U\not\equiv 0$. If $\{W_{m}\}$ is as in \eq{19}, then 
by \eq{4}, \eq{13}, \eq{18} and for large $m$
\begin{eqnarray}
&\di \frac{c}{2}\leq I(U_{m})-\frac{1}{2}\langle(-\Delta u_{m}^{1}+a(x)u_{m}^{1}-w_{m}^{1},
-\Delta u_{m}^{2}+b(x)u_{m}^{2}-w_{m}^{2}),U_{m}\rangle=\nonumber \\
&\di \frac{1}{2}\langle W_{m}, U_{m}\rangle -\int\limits_{\rn}F(x,U_{m})\, dx\leq\label{20}\\
&\di \frac{1}{2}\left(\int\limits_{u_{1}\geq 0}u_{1}(x)\uf(x,U)dx+
\int\limits_{u_{1}\leq 0}u_{1}(x)\of(x,U)dx+\int\limits_{u_{2}\geq 0}u_{1}(x)\ug(x,U)dx\
+\int\limits_{u_{2}\leq 0}u_{2}(x)\og(x,U)dx\right)\,.\nonumber
\end{eqnarray}
Now, taking into account the definition of $\of,\,\uf\,,\og,\,\ug$ we deduce that
$\of,\,\uf,\,\og,\,\ug$ verify \eq{15}, too.
So by \eq{20} we obtain
$$\frac{c}{2}\leq\int\limits_{\rn}(\eps|U_{m}|^{2}+A_{\eps}|u_{m}|^{p+1})=
\eps\|U_{m}\|_{L^{2}}^{2}+A_{\eps}\|U_{m}\|_{L^{p+1}}^{p+1}\,.$$
So, $\{U_{m}\}$ does not converge strongly to $0$ in $L^{p+1}(\rn;\r2)$.
From now on, a standard argument implies that $U\not\equiv 0$, which concludes our proof.\qed

\end{document}